\documentclass[12pt]{article}
\usepackage{amsmath}
\usepackage{amssymb}
\usepackage{amsthm}

\usepackage[margin=2.0cm]{geometry}
\usepackage{hyperref}
\usepackage{multirow}

\usepackage{color}
\usepackage{enumerate}
\usepackage{tikz,lmodern}
\usepackage{amsfonts}
\usetikzlibrary{arrows, decorations.markings}

\newtheorem{lem}{Lemma}[section]%
\newtheorem{theorem}[lem]{Theorem}%
\newtheorem{defi}[lem]{Definition}%
\newtheorem{cor}[lem]{Corollary}%
\newtheorem{prob}[lem]{Problem}%
\newtheorem{prop}[lem]{Proposition}%
\def\a{\alpha}
\def\b{\beta}
\def\g{\gamma}\def\d{\delta}

\def\s{\sigma}

\def\G{\Gamma}
\def\D{\bf D}

\def\lg{\langle}
\def\rg{\rangle}
\def\nd{\mathrel{\bigm|\kern-.7em/}}

\def\f{\noindent}

\def\Aut{\hbox{\rm Aut}}

\def\Cay{\hbox{\rm Cay}}

\def\Arc{\hbox{\rm Arc}}
\def\BCay{\hbox{\rm BCay}}

\def\demo{\f {\bf Proof.}\hskip10pt}

\def\mz{{\mathbb Z}}
\def\Z{{\mathbb Z}}
\def\D{\hbox{\rm D}}

\def\S{{\rm S}}

\def\G{\mathcal {G}}

\numberwithin{equation}{section}

\tikzstyle{vecArrow}=[semithick, decoration={markings,mark=at position 1 with {\arrow[semithick]{open triangle 60}}},
   double distance=1.4pt, shorten >= 5.5pt,
   preaction = {decorate},
   postaction = {draw,line width=1.4pt, white,shorten >= 4.5pt}]
\begin{document}
\title{On kernel isomorphisms of $m$-Cayley digraphs and finite $2$PCI-groups}
\author{\\Xing Zhang, Yan-Quan Feng\footnotemark, Jin-Xin Zhou, Fu-Gang Yin \\
{\small\em School of mathematics and statistics, Beijing Jiaotong University, Beijing, 100044, P.R. China}\\
}

\renewcommand{\thefootnote}{\fnsymbol{footnote}}
\footnotetext[1]{Corresponding author.
E-mails:
xingzh@bjtu.edu.cn (X. Zhang),
yqfeng@bjtu.edu.cn (Y.-Q. Feng),
jxzhou@bjtu.edu.cn (J.-X. Zhou),
fgyin@bjtu.edu.cn (F.-G. Yin)}

\date{}
\maketitle

\begin{abstract}
The isomorphism problem for digraphs is a fundamental problem in graph theory. In this paper, we consider this problem for $m$-Cayley digraphs which are generalization of Cayley digraphs. Let $m$ be a positive integer. A digraph admitting a group $G$ of automorphisms acting semiregularly on the vertices with exactly $m$ orbits is called an {\em $m$-Cayley digraph} of $G$. In our previous paper, we developed a theory for $m$-Cayley isomorphisms of $m$-Cayley digraphs, and classified finite $m$CI-groups for each $m\geq 2$, and finite $m$PCI-groups for each $m\geq 4$. The next natural step is to classify finite $m$PCI-groups for $m=2$ or $3$. Note that BCI-groups form an important subclass of the $2$PCI-groups, which were introduced in 2008 by Xu {\em et al.}. Despite much effort having been made on the study of BCI-groups, the problem of classifying finite BCI-groups is still widely open.

In this paper, we prove that every finite $2$PCI-group is solvable, and its Sylow $3$-subgroup is isomorphic to $\Z_3, \Z_3\times\Z_3$ or $\Z_9$, and Sylow $p$-subgroup with $p\not=3$ is either elementary abelian, or isomorphic to $\Z_4$ or $Q_8$.
We also introduce the kernel isomorphisms of $m$-Cayley digraphs, and establish some useful theory for studying this kind of isomorphisms. 
Using the results of kernel isomorphisms of $m$-Cayley digraphs together with the results on $2$PCI-groups, we give a proper description of finite BCI-groups, and in particular, we obtain a complete classification of finite non-abelian BCI-groups.
\bigskip

\noindent {\bf Keywords:}  $m$-Cayley digraph, $m$-PCayley digraph, Cayley isomorphism, $m$-Cayley isomorphism, semiregular group. \\
{\bf 2010 Mathematics Subject Classification:} 20B25, 05C25.
\end{abstract}

\section{Introduction}

Throughout this paper, all (di)graphs are finite and simple, and all groups are finite. For a digraph $\Gamma$, $\Aut(\Gamma)$ denotes the automorphism group of $\Gamma$, and $V(\Gamma)$ and $\Arc(\Gamma)$ denote the vertex set and arc set of $\Gamma$ respectively, while $E(\Gamma)$ denotes the edge set of $\Gamma$ when $\Gamma$ is a graph. For convenience, we view a graph $\Gamma$ with vertex set $V(\Gamma)$ and edge set $E(\Gamma)$ as a digraph with vertex set $V(\Gamma)$ and arc set $\{(u, v)\ :\ \{u,v\}\in E(\Gamma)\}$, and this makes no confusion in this paper. Thus, a graph is a special case of a digraph.

\subsection{$m$-Cayley digraphs and $m$-PCayley digraphs}
Let $G$ be a permutation group on a set $\Omega$. For $\a\in\Omega$, the stabilizer of $\a$ in $G$ is the subgroup $G_\a=\{g\in G\mid \a^g=\a\}$. We say that $G$ is {\em semiregular} on $\Omega$ if $G_\a=1$ for every $\a\in\Omega$. For a positive integer $m$, a digraph $\Gamma$ is called an {\em $m$-Cayley digraph of a group $G$} if $\Aut(\Gamma)$ has a semiregular subgroup which is isomorphic to $G$ and is semiregular on $V(\Gamma)$ with $m$ orbits. Note that every $m$-Cayley digraph $\Gamma$ of a group $G$ admits the following concrete realization. Let $S_{i,j} (1\leq i,j\leq m)$ be $m^2$ subsets of $G$. Define the digraph $\Gamma=\Cay(G, S_{i,j}: 1\leq i,j\leq m)$ to have vertex set
\[V(\Gamma)=\bigcup_{1\leq i\leq m}  G_i, \text{ where } G_i=\{x_i :  x\in G\},\]
and arc set
\[\mathrm{Arc} (\Gamma)=\bigcup_{1\leq i,j\leq m} \{(x_i, (sx)_j) :  s\in S_{i,j},x\in G\}.\]
We say that $\Gamma=\Cay(G, S_{i,j}: 1\leq i,j\leq m)$ is an {\em $m$-Cayley digraph of $G$ with respect to $S_{i,j} (1\leq i,j\leq m)$}. It can be easily seen that an $m$-Cayley digraph $\Gamma$ with the symbol $\{S_{i, j}: 1\leq i,j\leq m\}$ is undirected if and only if $S_{j, i}=S_{i, j}^{-1}=\{x^{-1}\mid x\in S_{i,j}\}$ for all $1\leq i,j\leq m$. For each $g\in G$, let $R(g)$ be the permutation on $V(\Gamma)$ defined by:
\[R(g): x_i \mapsto (xg)_i, \text{ for all } x  \in G \text{ and } 1\leq i \leq m.\]
Let $R(G)=\{R(g) : g \in G \}$. Then $R(G)\cong G$ and $R(G)$ acts semiregularly on $V(\Gamma)$ with $G_i (1\leq i\leq m)$ as its $m$ orbits.

By definition, when $m>1$ and $S_{i,i}=\emptyset$ for all $1\leq i \leq m$, the $m$-Cayley digraph $\Gamma=\Cay(G, S_{i,j}: 1\leq i,j\leq m)$ is a multipartite graph with $m$ parts $G_i (1\leq i\leq m)$, and in this case, we say that $\Gamma=\Cay(G, S_{i,j}: 1\leq i,j\leq m)$ is an \emph{$m$-PCayley digraph}.

\subsection{$m$-Cayley isomorphisms of $m$-Cayley (di)graphs}
In graph theory, one of the fundamental problems is to decide whether or not two given graphs are isomorphic. This paper is an attempt to investigate this problem for $m$-Cayley digraphs of the same group. Let $\Gamma=\Cay(G, S_{i,j}: 1\leq i,j\leq m)$ be an $m$-Cayley digraph of $G$. Let $N_{\mathrm{S}_{V(\Gamma)}}(R(G))$ be the normalizer of $R(G)$ in the symmetric group $\mathrm{S}_{V(\Gamma)}$ on $V(\Gamma)$. Take $\a\in N_{\mathrm{S}_{V(\Gamma)}}(R(G))$. It is easy to see that $\a$ is an isomorphism from $\Gamma$ to $\Gamma^\a=\Cay(G, T_{i,j}: 1\leq i,j\leq m)$, where $T_{i,j}=S_{i,j}^\a$ for $1\leq i,j\leq m$. We say that such an isomorphism $\a$ is an {\em $m$-Cayley isomorphism}. It is natural to investigate the conditions under which two $m$-Cayley (di)graphs of the same group are isomorphic if and only if there is an $m$-Cayley isomorphism between them.

\begin{defi}\label{def:mCI} {\rm Let $G$ be a group.
\begin{itemize}
\item An $m$-Cayley (di)graph $\Gamma$ of a group $G$ is said to be \emph{$m$CI} ($m$CI stands for $m$-Cayley isomorphism) if, for any $m$-Cayley (di)graph $\Sigma$ of $G$ isomorphic to $\Gamma$, there exists some $\a\in  N_{\mathrm{S}_{V(\Gamma)}}(R(G))$ such that $\Gamma^\a=\Sigma$.
\item If every $m$-Cayley digraph ($m$-Cayley graph, resp.) of $G$ is $m$CI, then $G$ is called an $m$DCI-group ($m$CI-group, resp.).
\end{itemize}}
\end{defi}

In case $m=1$, $1$-Cayley (di)graph is just the Cayley (di)graph, and $1$CI Cayley (di)graph is just the CI Cayley (di)graph. Furthermore, the $1$DCI-group and $1$CI-group are just DCI-group and CI-group, respectively. There are two long standing open questions on Cayley (di)graphs: (1)\ Which Cayley (di)graphs of a group $G$ are CI? (2)\ Which groups are DCI-groups or CI-groups?
There has been quite a lot of research on these two questions, for which we refer the reader to the nice survey paper~\cite{Li} and some recent papers~\cite{DE1,DT,Dobson2018-2,Dobson2018,FK,Kov,C.H.Li8,Mu3,Muzychuk2021, Ryabov2020,Ryabov2021,Xie2}.

Note that it is possible that an $m$-PCayley digraph of a group $G$ may be isomorphic to an $m$-Cayley digraph of a group $G$ which is not $m$-PCayley. So an interesting problem is to consider the isomorphisms between two $m$-PCayley (di)graphs of $G$. This motivates us to introduce the following concepts.

\begin{defi}\label{def:mPCI} {\rm Let $G$ be a group and let $m>1$.
\begin{itemize}
\item An $m$-PCayley (di)graph $\Gamma$ of $G$ is said to be \emph{mPCI} if, for any $m$-PCayley (di)graph $\Sigma$ of $G$ such that there is an isomorphism from $\Gamma$ to $\Sigma$ keeping $\{G_1, G_2, \cdots, G_m\}$ invariant, there is some $n \in  N_{\mathrm{S}_{V(\Gamma)}}(R(G))$ such that $\Gamma^n=\Sigma$.
\item If every $m$-PCayley digraph ($m$-PCayley graph, resp.) of $G$ is $m$PCI, then $G$ is called an $m$PDCI-group ($m$PCI-group, resp.).
\end{itemize}}
\end{defi}

The $2$CI-(di)graphs and $2$PCI-(di)graphs have also been received considerable attention (see, for example, \cite{Arez,Dobson3}). In our previous paper~\cite{ZFYZ}, we developed some general theory for $m$CI-(di)graphs and $m$PCI-(di)graphs. With these, we obtained, somewhat to our surprise, a complete classification of finite $m$CI- and $m$DCI-groups for each $m\geq 2$, and a complete classification of finite $m$PCI- and $m$PDCI-groups for each $m\geq 4$. This suggests that the essential difficult case in the study of $m$CI- and $m$DCI-groups is $m=1$, namely, the ordinary CI- and DCI-groups. For $m$PCI- and $m$PDCI-groups, the remaining cases are $m=2$ and $m=3$. In this paper, we shall deal with $2$PCI-groups. The following is our first main theorem.

\begin{theorem}\label{Sylow2pci}
Let $G$ be a $2$PCI-group. Then $G$ is solvable, and every Sylow $3$-subgroup of $G$ is isomorphic to $\Z_3, \Z_3\times\Z_3$ or $\Z_9$, and every Sylow $p$-subgroup of $G$ with $p\neq 3$ is either elementary abelian, or is isomorphic to $\Z_4$ or $Q_8$.
\end{theorem}

This theorem implies that every finite $2$PCI-group is solvable and has very restrictive subgroup structure. So we would like to propose the following problem.

\begin{prob}\label{Openproblem}
Classify finite $2$PCI-groups.
\end{prob}

\subsection{Kernel $m$-Cayley isomorphisms of $m$-Cayley (di)graphs}

Xu et al.~in \cite{Xu} initiated the study of BCI-groups, which form a special class of $2$PCI-groups defined below. Let $\Gamma$ be a $2$-PCayley graph of a group $G$. Let $N=N_{\mathrm{S}_{V(\Gamma)}}(R(G))$ and let $K$ be the kernel of $N$ acting on $\{G_1, G_2\}$. Define $K_{(1_G)_1}$ as the stabilizer of the vertex $(1_G)_1$ in $K$. Then $\Gamma$ is said to be {\em BCI} if, for any $2$-PCayley graph $\Sigma$ of $G$ isomorphic to $\Gamma$, there exists some $\a\in  K_{(1_G)_1}$ such that $\Gamma^\a=\Sigma$. As $R(G)\leq\Aut(\Gamma)$, according to Proposition~\ref{c-kernel}, it is not difficult to see that $\Gamma$ is BCI if and only if for any $2$-PCayley graph $\Sigma$ of $G$ isomorphic to $\Gamma$, there exists some $\a\in  K$ such that $\Gamma^\a=\Sigma$.

In general, let $\Gamma$ be an $m$-PCayley (di)graph of a group $G$ with $m\geq2$. Let $K$ be the kernel of $N_{\mathrm{S}_{V(\Gamma)}}(R(G))$ acting on $\{G_1, G_2, \dots, G_m\}$. For any $\b\in K$, it is easy to see that $\Gamma^\b$ is also an $m$-PCayley (di)graph of $G$. Clearly, $\b$ is an isomorphism from $\Gamma$ to $\Gamma^\b$. Such an isomorphism $\b$ will be called a {\em kernel $m$-Cayley isomorphism}. Motivated by this, we introduce the following concepts which generalize the BCI-graphs.

\begin{defi}\label{def:mPCI}
{\rm Let $\Gamma$ be an $m$-Cayley (di)graph of a group $G$. Let $N=N_{\mathrm{S}_{V(\Gamma)}}(R(G))$ and let $K$ be the kernel of $N$ acting on the orbits set $\{G_1, G_2, \dots, G_m\}$ of $R(G)$ on $V(\Gamma)$.
\begin{itemize}
\item $\Gamma$ is said to be \emph{K$m$CI} (stands for kernel $m$-Cayley isomorphism) if, for any $m$-Cayley (di)graph $\Sigma$ of $G$ isomorphic to $\Gamma$, there exists $k\in K$ such that $\Gamma^k=\Sigma$.
\item If every $m$-Cayley digraph ($m$-Cayley graph, resp.) of $G$ is K$m$CI, then $G$ is called a K$m$DCI-group (K$m$CI-group, resp.).
\end{itemize}
Assume that $\Gamma$ is an $m$-PCayley (di)graph of a group $G$. Then
\begin{itemize}
\item $\Gamma$ is said to be \emph{K$m$PCI} if, for any $m$-PCayley (di)graph $\Sigma$ of $G$ such that there is an isomorphism from $\Gamma$ to $\Sigma$ keeping $\{G_1,G_2,\cdots,G_m\}$ invariant, there is some $k\in K$ such that $\Gamma^k=\Sigma$.
\item If every $m$-Cayley digraph ($m$-Cayley graph, resp.) of $G$ is K$m$PCI, then $G$ is called an K$m$PDCI-group (K$m$PCI-group, resp.).
\end{itemize}}
\end{defi}

By definition, every K$m$PCI-(di)graph is also an $m$PCI-(di)graph. In case $m=2$, the K2PCI-graphs (namely the BCI-graphs) have been extensively studied in the literature. For example, a criterion for a $2$-PCayley graph being K2PCI was given in~\cite[Theorem~C]{Arez} (see also~\cite[Lemma~2.21]{Ko2}). In~\cite{Are1}, it was proved that every K2PCI-group (namely, BCI-group) is solvable. In~\cite{Arez}, it was proved that every cyclic group of order $9$ or a prime order is K2PCI, and moreover, it was proved that every Sylow $p$-subgroup of a K2PCI-group is elementary abelian for $p>3$. There is also some work on K2PCI-graphs with valency at most $3$ (see, for example, \cite{Ji2,Ko1,Ko2}). For more work about K2PCI-graphs, we refer the reader to~\cite{Are2,Ji1,Ji3,Ji4,Ko4,Ko3}.

Motivated by the facts listed above, in this paper we shall investigate the kernel isomorphisms of $m$-Cayley (di)graphs over the same group. We first give a criterion for an $m$-Cayley (di)graph being K$m$CI as well as a criterion for an $m$-PCayley (di)graph being K$m$PCI. With these, we prove the following theorem which gives a classification of K$m$CI-groups and K$m$DCI-groups for all $m\geq2$, and a classification of K$m$PCI-groups and K$m$PDCI-groups for all $m\geq3$.

\begin{theorem}\label{mCI-groups}
Let $m$ be a positive integer and let $G$ be a finite group. Then
  \begin{itemize}
    \item[\rm (1)] For $m\geq 2$, $G$ is K$m$CI if and only if $m=2$ and $G=1$;
    \item[\rm (2)] For $m\geq 2$, $G$ cannot be K$m$DCI.
    \item[\rm (3)] For $m\geq 3$, $G$ cannot be K$m$PCI and K$m$PDCI.
 \end{itemize}
\end{theorem}

This theorem reveals that the main difficult in the study of K$m$PCI-groups and K$m$PDCI-groups is the case $m=2$. Our third main result gives a characterization of abelian K$2$PCI-groups, and in particular, a complete classification is given of non-abelian K$2$PCI-groups.

\begin{theorem}\label{BCIEA}
Let $G$ be a group. Then the following hold.
\begin{enumerate}[{\rm (1)}]
\item\ If $G$ is abelian, then $G$ is K$2$PCI (namely, BCI) if and only if every Sylow $3$-subgroup of $G$ is isomorphic to $\Z_3, \Z_3\times\Z_3$ or $\Z_9$, every Sylow $p$-subgroup of $G$ with $p\neq 3$ is either elementary abelian, or is isomorphic to $\Z_4$, and for any $S\subseteq G$,
    $\Cay(\langle S\rangle, S,\emptyset,\emptyset,S^{-1})$ is K$2$PCI.

\item\ If $G$ is non-abelian, then $G$ is K$2$PCI (namely, BCI) if and only if $G$ is isomorphic to one of the following groups: $D_6, Q_8$ or $D_{10}$.
\end{enumerate}
\end{theorem}

\medskip
Based on this theorem, to classify finite K$2$PCI (namely, BCI)-groups, it suffices to classify all abelian K$2$PCI-groups. We leave this as an open problem.

\begin{prob}\label{Openproblem}
Classify abelian K$2$PCI-groups.
\end{prob}

\medskip
The paper is organised as follows. In Section~\ref{sec2}, we prove the criterions for K$m$CI- and K$m$PCI-(di)graphs, and then using this, we prove Theorem~\ref{mCI-groups}. In Section~\ref{sec3}, we investigate the $2$-PCI groups, and give some basic properties as well as some sufficient and necessary conditions for $2$PCI-groups, and then using these, we prove Theorem~\ref{Sylow2pci} and Theorem~\ref{BCIEA}.

To end this section we fix some notation used in this paper. The notation for finite groups in this paper is standard; see \cite{Atlas} for example. In particular, denote by $\mz_{m}$ the additive group of integer numbers modulo $m$, by $\D_{2m}$ the dihedral group of order $2m$, and by $\S_m$ the symmetric group of degree $m$. Also, we denote by $\S_{\Omega}$ the symmetric group on the set $\Omega$.
For a prime $p$,  we use $\mz_p^m$ to denote the elementary abelian group of order $p^m$. For two groups $A$ and $B$, $A\times B$ stands for the direct product of $A$ and $B$, and $A:B$ or $A\rtimes B$ for a split extension or a semi-direct product of $A$ by $B$.

\section{Kernel isomorphisms of $m$-Cayley (di)graphs}\label{sec2}

\subsection{Preliminaries}\label{notation}

Throughout this section, we fix the following notation.
\begin{itemize}
  \item $m$:\ a positive integer
  \item $G$:\ a finite group
  \item $\Gamma=\Cay(G,S_{i,j}: 1\leq i,j\leq m)$:\ the $m$-Cayley digraph of $G$ with respect to the subsets $S_{i,j}$ for $1\leq i,j\leq m$
  \item $\G=\{G_1, G_2, \dots, G_m\}$:\ the set of $R(G)$ orbits on $V(\Gamma)$
  \item $N=N_{\S_{V(\Gamma)}}(R(G))$:\ the normalizer of $R(G)$ in the symmetric group $\S_{V(\Gamma)}$
  \item $K$:\ the kernel of $N$ acting on the set of $R(G)$ orbits on $V(\Gamma)$
\end{itemize}

For $g\in G$ and $1\leq i\leq m$, the permutation $L_i(g)$ on $V(\Gamma)$ is define by
\[
 L_i(g):\ x_{i}\mapsto (g^{-1}x)_{i}, \text{ and } x_{j}\mapsto x_{j} \text{ for all } x\in G \text{ and } 1\leq j\leq m \text{ with } j\neq i.\]
Set $L_i=\{L_i(g): g\in G\} \text{ and } L=L_1L_2\cdots L_m$. Then $L_i$ is regular on $G_i$ and fixes $\bigcup\limits_{j\not=i}^{m}G_{j}$ pointwise. It is easy to see that for all $1\leq i,j\leq m$ with $i\neq j$, $L_i\cap L_j=1$ and $L_i$ and $L_j$ commute pointwise. Thus, $L = L_1\times L_2\times\cdots\times L_m$. Clearly, $|L_i|=|R(G)|=|G|$ and $|L|=|G|^m$.
Define two natural permutations on $V(\Gamma)$: for every $\sigma\in \S_{m}$ (the symmetric group on $\{1,2,\cdots,m\}$), $\alpha \in \mathrm{Aut}(G)$, $i\in \{1,2,\ldots,m\}$ and $x\in G$:
\[
\sigma: x_i \mapsto x_{i^\sigma},\ \ \ \ \ \ \
\alpha: x_i \mapsto (x^\alpha)_i.\]

The following proposition determines the groups $N,K$ and $K_{1_r}$.

\begin{prop}[{\cite[Corollaries 2.4--2.6]{ZFYZ}}]\label{c-kernel}
$(1)$\ The normalizer $N=(L_1\times\cdots\times L_m)\rtimes (\S_m\times \Aut(G))=((L_1\times\cdots \times L_{r-1}\times R(G)\times L_{r+1}\times\cdots\times L_m)\rtimes \Aut(G))\rtimes\S_m$ for every $1\leq r\leq m$.

$(2)$\ The kernel $K=L\rtimes\Aut(G)=(L_1\times\cdots \times L_{r-1}\times R(G)\times L_{r+1}\times\cdots\times L_m)\rtimes \Aut(G)$  for every $1\leq r\leq m$.

$(3)$\ The stabilizer $K_{1_r}=(L_1\times \cdots \times L_{r-1}\times L_{r+1}\times\cdots\times L_m)\rtimes \Aut(G)$ for every $1\leq r\leq m$.
\end{prop}

Using the above proposition, we can obtain the following proposition which is very useful in the study of $\mathcal{P}$-(di)graphs, where $\mathcal{P}\in \{mPCI, KmCI, KmPCI\}$.

\begin{prop}\label{mCImPCI}
Let $1\leq r\leq m$. Then the following hold.
\begin{enumerate}[{\rm (1)}]
  \item $\Gamma=\Cay(G,S_{i,j}: 1\leq i,j\leq m)$ is K$m$CI if and only if whenever $\Gamma\cong\Cay(G, T_{i,j}: 1\leq i,j\leq m)$, there exist $\alpha\in\Aut(G)$ and $g_1,\ldots, g_m\in G$ with $g_{r}=1$ such that $T_{i,j}=g_jS^{\alpha}_{i,j}g_{i}^{-1}$.

  \item Assume that $\Gamma=\Cay(G,S_{i,j}: 1\leq i,j\leq m)$ is an $m$-PCayley digraph. Then the following hold.
  \begin{enumerate}[{\rm (i)}]
    \item $\Gamma$ is $m$PCI if and only if for any $m$-PCayley (di)graph $\Sigma=\Cay(G, T_{i,j}: 1\leq i,j\leq m)$ such that there is an isomorphism from $\Gamma$ to $\Sigma$ keeping $\G$ invariant, there exist $\alpha\in\Aut(G)$, $\s\in\S_m$ and $g_1,\ldots, g_m\in G$ with $g_{r}=1$ such that $T_{i^\s,j^\s}=g_jS^{\alpha}_{i,j}g_{i}^{-1}$.
    \item $\Gamma$ is K$m$PCI if and only if for any $m$-PCayley (di)graph $\Sigma=\Cay(G, T_{i,j}: 1\leq i,j\leq m)$ such that there is an isomorphism from $\Gamma$ to $\Sigma$ keeping $\G$ invariant, there exist $\alpha\in\Aut(G)$ and $g_1,\ldots, g_m\in G$ with $g_{r}=1$ such that $T_{i,j}=g_jS^{\alpha}_{i,j}g_{i}^{-1}$.
  \end{enumerate}
\end{enumerate}
\end{prop}

\demo Take $n\in N$ and $1\leq r\leq m$. By Proposition~\ref{c-kernel}, there are $g,g_1,g_2,\cdots,g_m\in G$ with $g_r=1$, $\s\in\S_m$ and $\a\in \Aut(G)$, such that
\[n=L_1(g_1)\cdots L_{r-1}(g_{r-1})R(g)L_{r+1}(g_{r+1})\cdots L_m(g_m)\a\s.\]
Note that $n\in K$ if and only if $\s=1$.

For any arc $(x_i, (s_{i,j}x)_j)\in \Arc(\Gamma)$, where $x\in G$, $1\leq i, j\leq m$ and $s_{i,j}\in S_{i,j}$, we have
\begin{equation*}\label{EqninN}
(x_i,(s_{i,j}x)_j)^n=((g_i^{-1}x)^\a_{i^\s},(g_j^{-1}s_{i,j}x)^\a_{j^\s})=((g_i^{-1}x)^\a_{i^\s},((g_j^{-1}s_{i,j}g_i)^\a (g_i^{-1}x)^\a)_{j^\s}).
\end{equation*}
Then
\begin{equation*}\label{Eqmkdci}
\Gamma^n=\Cay(G, T_{i,j}: 1\leq i,j\leq m), \mbox{ where } T_{i^\s, j^\s}=(g_j^{-1}S_{i,j}g_i)^\a.
\end{equation*}
By the definition of $\mathcal{P}$-(di)graphs with $\mathcal{P}\in\{mPCI, KmCI, KmPCI\}$, we immediately obtain this proposition.\hfill\qed

\subsection{Criterions for K$m$(P)CI $m$-Cayley (di)graphs\label{sec3}}

We first give a criterion for an $m$-Cayley graph being K$m$CI.

\begin{theorem}\label{BabaiSimilar}
Let $\Gamma=\Cay(G,S_{i,j}: 1\leq i,j\leq m)$ be an
$m$-Cayley (di)graph of a group $G$. Then $\Gamma$ is K$m$CI if and only if the normalizer $N_{\mathrm{Aut}(\Gamma)}(R(G))$ induces the symmetric group  on $\{G_1,G_2,\cdots,G_m\}$, and every semiregular subgroup of $\Aut(\Gamma)$ isomorphic to $G$ is conjugate to $R(G)$ in $\Aut(\Gamma)$.
\end{theorem}

\demo Recall that $V(\Gamma)=\bigcup_{i=1}^m G_i$ and $\G=\{G_1,G_2,\cdots,G_m\}$. By Proposition~\ref{c-kernel} $N=(L_1\times L_2\times\cdots\times L_m)\rtimes(\Aut(G)\times \S_m)$. Then the symmetric group $\S_\G$ on $\G$ is isomorphic to $\S_m$. Set $A=\Aut(\Gamma)$.

To prove the necessity, assume that $\Gamma$ is K$m$CI and let $H\cong R(G)$ be a semiregular subgroup of $A$ with $m$ orbits on $V(\Gamma)$.
By~\cite[Theorem 2.8]{ZFYZ}, $H$ and $R(G)$ are conjugate in $\S_{V(\Gamma)}$, and then there exists $\s\in\S_{V(\Gamma)}$ such that $R(G)=H^\s$.
Note that $\s$ is an isomorphism from $\Gamma$ to $\Gamma^\s$. Since $H$ is a semiregular subgroup of $A=\Aut(\Gamma)$, $H^\s$ is a semiregular subgroup of $\Aut(\Gamma^\s)$, and since $R(G)=H^\s$, we have $R(G)\leq \Aut(\Gamma^\s)$. This implies that $\Gamma^\s$ is also an $m$-Cayley (di)graph of $G$. Since $\Gamma$ is K$m$CI, there exists a $k\in K$ such that $\Gamma^k=\Gamma^\s$. Thus, $k\s^{-1}\in\Aut(\Gamma)=A$ and $R(G)^{k\s^{-1}}=R(G)^{\s^{-1}}=H$, that is, $H$ and $R(G)$ are conjugate in $A$. We are left to show that $N_A(R(G))$ induces the symmetric group $\S_{\G}\cong\S_m$. To do so, for any $\s\in\S_m$, we shall seek a $\b\in N_A(R(G))$ such that $G_i^\b=G_{i^\s}$ for all $1\leq i\leq m$.
Note that $N=(L_1\times L_2\times\cdots\times L_m)\rtimes(\Aut(G)\times \S_m)$ induces the symmetric group $\S_\G$ on $\G$. Then there is an $n\in N$ such that $G_i^{n^{-1}}=G_{i^\s}$ for all $1\leq i\leq m$. Clearly, $n$ induces an isomorphism from $\Gamma$ to $\Gamma^n$, and hence $R(G)=R(G)^n\leq \Aut(\Gamma^n)$, implying that $\Gamma^n$ is an $m$-Cayley (di)graph of $G$. Since $\Gamma$ is K$m$CI, there is a $k\in K=(L_1\times L_2\times\cdots\times L_m)\rtimes \Aut(G)$ such that $\Gamma^k=\Gamma^n$. Thus, $kn^{-1}\in A=\Aut(\Gamma)$, and since $kn^{-1}\in N$, we have $kn^{-1}\in N_A(R(G))$. Take $\b=kn^{-1}$. Then $G_i^\b=G_i^{kn^{-1}}=G_i^{n^{-1}}=G_{i^\s}$ for all $1\leq i\leq m$, as required.

We now prove the sufficiency. Assume that the normalizer $N$ induces the symmetric group $\S_\G$ on $\G$ and every semiregular subgroup of $A$ isomorphic to $G$ is conjugate to $R(G)$ in $A$. Let $\Sigma$ be an $m$-Cayley (di)graph of $G$ that is isomorphic to $\Gamma$. To finish the proof, it suffices to show that there is $k\in K$ such that $\Gamma^k=\Sigma$. Since $\Gamma\cong \Sigma$, we may assume that $\g$ is an isomorphism from $\Gamma$ to $\Sigma$. Since $\Sigma$ is an $m$-Cayley (di)graph of $G$, we have $R(G)\leq \Aut(\Sigma)$, and since $\Gamma^\g=\Sigma$, we have $R(G)^{\g^{-1}}\leq A$. Clearly, $R(G)$ and $R(G)^{\g^{-1}}$ are isomorphic and semiregular on $V(\Gamma)$. By assumption, there is $\a\in A$ such that $R(G)^\a=R(G)^{\g^{-1}}$. Thus, $\a\g\in N$. By Corollary~\ref{c-kernel}, $K=L\rtimes\Aut(G)$ is the kernel of $N=L\rtimes(\Aut(G)\times \S_m)$ acting on $\G$, where $L=L_1\times L_2\times\cdots\times L_m$.  Thus, $N/K\cong \S_m$. Since $N_A(R(G))\leq N$ induces the symmetric group $\S_\G\cong\S_m$, we have $KN_A(R(G))/K\cong \S_m$. It follows that $N=N_A(R(G))K$. Thus, $\a\g=ak$ for some $a\in N_A(R(G))$ and $k\in K$.  Since $a,\a\in A$, we have $\Gamma^a=\Gamma$ and $\Gamma^\a=\Gamma$, and hence   $\Gamma^k=\Gamma^{ak}=\Gamma^{\a\g}=\Gamma^\g=\Sigma$, completing the proof. \hfill\qed

\medskip
Next we give a criterion for an $m$-PCayley (di)graph being K$m$PCI.

\begin{theorem}\label{PCI-BabaiSimilar}
Let $\Gamma=\Cay(G,S_{i,j}: 1\leq i,j\leq m)$ be an
$m$-PCayley (di)graph. Then $\Gamma$ is K$m$PCI if and only if the normalizer $N$ induces the symmetric group on $\G$, and every semiregular subgroup of $\Aut(\Gamma)$, which is isomorphic to $G$ and has the same orbit set $\G$ with $R(G)$, is conjugate to $R(G)$ in $\Aut(\Gamma)$.
\end{theorem}

\demo  Let $A=\Aut(\Gamma)$. To prove the necessity, assume that $\Gamma$ is K$m$PCI and let $H$ be a semiregular subgroup of $\Aut(\Gamma)$, is isomorphic to $G$ and has the same orbit set $\G$ with $R(G)$. According to~\cite[Theorem 2.8]{ZFYZ}, there exists a $\s\in\S_{V(\Gamma)}$ such that $R(G)=H^\s$. Since $H$ has the orbit set $\G$ on $V(\Gamma)$, $H^\s$ has the orbit set $\{G_1^\s, G_2^\s, \cdots, G_m^\s\}$ on $V(\Gamma)$. This implies that $\G=\{G_1^\s, G_2^\s, \cdots, G_m^\s\}$. So $\Gamma^\s$ is also an $m$-PCayley (di)graph of $G$. Since $\Gamma$ is K$m$PCI, there is a $k\in K$ such that $\Gamma^k=\Gamma^\s$. It follows $k\s^{-1}\in A$ and $R(G)^{k\s^{-1}}=R(G)^{\s^{-1}}=H$, that is $H$ and $R(G)$ are conjugate in $A$.  We are only left to show that $N_A(R(G))$ induces the symmetric group $\S_{\G}$. To do this, for $\s\in\S_m$, we shall prove that there is a $\b\in N_A(R(G))$ such that $G_i^\b=G_{i^\s}$ for all $1\leq i\leq m$. Since $N$ induces the symmetric group $\S_\G$ on $\G$, there is $n\in N$ such that $G_i^{n^{-1}}=G_{i^\s}$ for all $1\leq i\leq m$. Clearly, $n$ keeps $\G$ invariant, so $\Gamma^n$ is also an $m$-PCayley (di)graph of $G$. Since $\Gamma$ is K$m$CI, there is a $k\in K=(L_1\times L_2\times\cdots\times L_m)\rtimes \Aut(G)$ such that $\Gamma^k=\Gamma^n$. Write $\b=kn^{-1}$. Then $\b\in A$, and so $\b\in N$. Furthermore, $G_i^\b=G_i^{kn^{-1}}=G_i^{n^{-1}}=G_{i^\s}$, as required.

We now prove the sufficiency. Assume that $N$ induces the symmetric group $\S_\G$ on $\G$ and that every semiregular subgroup of $A$, which is isomorphic to $G$ and has orbit set $\G$ on $V(\Gamma)$, is conjugate to $R(G)$ in $\Aut(\Gamma)$. Let $\Sigma$ be an $m$-PCayley (di)graph of $G$ and assume that $\g$ is an isomorphism $\Gamma$ to $\Sigma$ keeping $\G$ invariant. To finish the proof, it suffices to show that there is a $k\in K$ such that $\Gamma^k=\Sigma$. Since $\g$ is an isomorphism $\Gamma$ to $\Sigma$ keeping $\G$ invariant, $R(G)$ and $R(G)^{\g^{-1}}$ are two isomorphic semiregular subgroups of $\Aut(\Gamma)$, both of which have the same orbit set $\G$ on $V(\Gamma)$. By assumption, there is $\a\in A$ such that $R(G)^\a=R(G)^{\g^{-1}}$, implying $\a\g\in N$. Since $N_A(R(G))\leq N$ induces the symmetric group $\S_\G\cong\S_m$, we have $KN_A(R(G))/K\cong \S_m$ and so $N=N_A(R(G))K$. Thus, $\a\g=ak$ for some $a\in N_A(R(G))$ and $k\in K$.  Since $a,\a\in A$, we have $\Gamma^a=\Gamma$ and $\Gamma^\a=\Gamma$, and hence $\Gamma^k=\Gamma^{ak}=\Gamma^{\a\g}=\Gamma^\g=\Sigma$, as required.  \hfill\qed

\medskip
\f{\bf Remark.}\ The case $m=2$ in Theorem~\ref{PCI-BabaiSimilar} has been proved in \cite[Theorem~C]{Arez} (see also \cite[Lemma 2.6]{Ko1} and \cite[Lemma 2.2]{Ko2}).
\medskip

\subsection{K$m$CI-, K$m$DCI-, K$m$PCI- and K$m$PDCI-groups}

In this section we shall give a classification of K$m$CI- and K$m$DCI-groups for every $m\geq 2$ as well as a classification of K$m$PCI- and K$m$PDCI-groups for each $m\geq3$. Before doing this, we first prove a lemma.

\begin{lem}\label{mtom-1}
Let $m\geq 2$ be an integer and let $G$ be a finite group. Let
\[\mathcal{P}_m\in\{mDCI, mCI, KmDCI, KmCI, mPDCI, mPCI, KmPDCI, KmPCI\}.\]
If $G$ is $\mathcal{P}_m$-group, then $G$ is also $\mathcal{P}_{m-1}$-group.
\end{lem}

\demo If $\mathcal{P}_m\in\{mDCI,mCI,mPDCI,mPCI\}$, the result follows from~\cite[Theorem 5.3]{ZFYZ}. We only need to deal with the cases where $\mathcal{P}_m\in\{KmDCI,KmCI,KmPDCI,KmPCI\}$.

Let $V_m=G_1\cup G_{2}\cup\dots\cup G_{m}$ and $V_{m-1}=G_1\cup G_{2}\cup\dots\cup G_{m-1}$. Let $K_{m-1}$ denote the kernel of $N_{\S_{V_{m-1}}}(R(G))$ acting on $\{G_1,G_2,\cdots,G_{m-1}\}$.

Assume that $G$ is K$m$DCI. To prove that $G$ is K$(m-1)$DCI, let $\Gamma=\Cay(G,S_{i,j}: 1\leq i,j\leq m-1)$ and $\Sigma=\Cay(G,T_{i,j}: 1\leq i,j\leq m-1)$ be two isomorphic $(m-1)$-Cayley digraphs of $G$, and we only need to prove that there is $k\in K_{m-1}$ such that $\Gamma^k=\Sigma$. Note that $S_{i,j}$ and $T_{i,j}$ are given for every $1\leq i,j\leq m-1$. We extend them to $1\leq i,j\leq m$ by setting $S_{k,m}=S_{m,k}=S_{m,m}=T_{k,m}=T_{m,k}=T_{m,m}=\emptyset$ for each
$1\leq k\leq m-1$. Let $\Gamma_1=\Cay(G,S_{i,j}: 1\leq i,j\leq m)$ and $\Sigma_1=\Cay(G,T_{i,j}: 1\leq i,j\leq m)$. Then $\Gamma_1=\Gamma\cup G_m$, $\Sigma_1= \Sigma\cup G_m$, and so $\Gamma_1$ and $\Sigma_1$ are isomorphic $m$-Cayley digraphs.  Since $G$ is a K$m$DCI-group, there is $k'\in K$ such that $\Gamma_1^{k'}=\Sigma_1$. Since $k'$ fixes every $G_i$ setwise for $1\leq i\leq m$, and hence $\Gamma^{k'}=\Sigma$. Let $k$ be the restriction of $k'$ on $V_{m-1}$. We have $k\in K_{m-1}$ and $\Gamma^k=\Sigma$, as required.

Similarly, one may prove that if $G$ is $\mathcal{P}_m$-group, then $G$ is also $\mathcal{P}_{m-1}$-group, where
\[\mathcal{P}_m\in\{KmCI, KmPDCI, KmPCI\}.\]
\hfill\qed
\medskip

\noindent{\bf Proof of Theorem~\ref{mCI-groups}:}\ For part (1), assume that $G$ is K$m$CI for $m\geq 2$. By definition, $G$ must be $m$CI, and by~\cite[Theorem 5.4(1)]{ZFYZ}, either $m=2$ and $G=1$ or $\mz_3$, or $m\geq 3$ and $G=1$.

Let $m=3$ and $G=1$. Take $S_{1,2}=S_{2,1}=\{1\}$ and $S_{i,j}=\emptyset$ otherwise. Then $\Aut(\Cay(G, S_{i,j}: 1\leq i,j\leq 3))\cong\mz_2$, which can not induce the symmetric group $\S_3$ on $\{G_1, G_2, G_3\}$. By Theorem~\ref{BabaiSimilar}, $G$ can not be K$3$CI. By Lemma~\ref{mtom-1}, $G$ cannot be K$m$CI for any $m\geq 3$.

Let $m=2$ and $G=\mathbb{Z}_3=\langle x\rangle$. Take $S_{1,1}=\{x,x^{-1}\}$ and $S_{1,2}=S_{2,1}=S_{2,2}=\emptyset$. Then $\Gamma:=\Cay(G,S_{i,j}: 1\leq i,j\leq 2)$ is a union of a triangle and three isolated vertices, and the triangle is just the subgraph induced by $G_1$. Thus, $\Aut(\Gamma)$ fixes $G_1$ and $G_2$ setwise. By Theorem~\ref{BabaiSimilar}, $\Gamma$ is not K$2$CI and so $G=\mathbb{Z}_3$ is not a K$2$CI-group.

Let $m=2$ and $G=1$. Then every $2$-Cayley graph $\Gamma$ of $G$ is either $2K_1$ or $K_2$. Thus $\Aut(\Gamma)=\mathbb{Z}_2$ and $\Aut(\Gamma)$ induces the symmetric group $\S_2$ on $\{G_1,G_2\}$. By Theorem~\ref{BabaiSimilar}, $\Gamma$ is K$2$CI, and so $G$ is K$2$CI. This completes the proof of part~(1).

For part (2), assume that $G$ is K$m$DCI for $m\geq 2$. By definition, $G$ must be K$m$CI, and by part~(1), $m=2$ and $G=1$. Take $S_{1,2}=\{1\}$ and $S_{1,1}=S_{2,2}=S_{2,1}=\emptyset$. Then $\Gamma:=\Cay(G,S_{i,j}: 1\leq i,j\leq 2)$ is an arc, and hence $\Aut(\Gamma)$ fixes $G_1$ and $G_2$ setwise. By Theorem~\ref{BabaiSimilar}, $G$ is not a K$2$DCI-group. This implies that
any finite group cannot be a K$m$DCI-group for $m\geq 2$. This completes the proof of part~(2).

For part (3), assume that $m\geq 3$. Take $S_{1,2}=S_{2,1}=\{1\}$, and $S_{i,j}=\emptyset$ for other $1\leq i,j\leq 3$. Then $\Gamma:=\Cay(G,S_{i,j}: 1\leq i,j\leq 3)\cong |G|K_2\cup |G|K_1$, and $\Aut(\Gamma)$ fixes $G_3$ setwise. By Theorem~\ref{PCI-BabaiSimilar}, $\Gamma$ is not K$3$PCI and so not K$m$PCI for every $m\geq 3$ by Lemma~\ref{mtom-1}.
By definition, a finite K$m$PDCI-group must be a K$m$PCI, and hence there is no finite K$m$PDCI-groups for every $m\geq 3$. This completes the proof of part~(3).
\hfill\qed

\section{Finite $2$PCI-groups\label{sec3}}

Let $G$ be a finite group and $\Gamma=\Cay(G, S_{1,1}, S_{1,2}, S_{2,1}, S_{2,2})$ be a $2$-PCayley graph of $G$. Then $S_{1,1}=S_{2,2}=\emptyset$ and $S_{1,2}=S_{2,1}^{-1}=\{x^{-1}\mid x\in S_{2,1}\}$. For convenience, throughout this section we shall follow \cite{Xu} to use the notation for $\Gamma$:
\[\BCay(G, S):=\Cay(G, S_{1,1}, S_{1,2}, S_{2,1}, S_{2,2}),\ {\rm where}\ S=S_{1,2}.\]

\subsection{Some basic properties of $2$-PCayley graphs}

We start with a lemma which gives a few basic properties of $2$-PCayley graphs.

\begin{lem}\label{basic}
Let $G$ be a finite group and $S, T\subseteq G$. Then the following hold.
\begin{itemize}
  \item[\rm (i)]$\ \BCay(G, S)$ is connected if and only if $G=\langle S S^{-1}\rangle$.
  \item[\rm (ii)]$\ \BCay(G, S) \cong \BCay(G, T)$ if and only if $\BCay(\langle S S^{-1}\rangle, S) \cong \BCay(\langle T T^{-1}\rangle, T)$.
  \item[\rm (iii)]\ Assume that $\BCay(G,S)$ is $2$PCI. If $\BCay(G, S)\cong \BCay(G, T)$, then $\BCay(G, T)$ is also $2$PCI.
  \item[\rm (iv)]\ If $\BCay(G,S)$ is $2$PCI, then up to graph isomorphism we may assume that $S$ contains the identity element of $G$.
\end{itemize}
\end{lem}

\demo Part (i) was proved in~\cite{DuXu}, and when $1 \in S$, it is easy to see $\langle S S^{-1}\rangle=\langle S\rangle$.
Part (ii) is obvious, see~\cite[Lemma 2.8]{Ji2}.

For part (iii), let $\Gamma=\BCay(G, S)$ and $\Sigma=\BCay(G,T)$. Assume that $\Gamma$ is $2$PCI. If $\Gamma\cong\Sigma$, then by \cite[Lemma~4.4]{Arez}, there always exists an isomorphism, say $\b$, from $\Gamma$ to $\Sigma$ such that $\b$ preserves $\{G_1, G_2\}$. Since $\Gamma$ is $2$PCI, there exists $n\in N_{\S_{V(\Gamma)}}(R(G))$ such that $\Gamma^n=\Sigma$. Take a $2$-PCayley graph $\Sigma'$ of $G$ such that there exists an isomorphism, say $\d$, from $\Sigma$ to $\Sigma'$ keeping $\{G_1, G_2\}$ invariant. Then we also have $\b\d$ is an isomorphism from $\Gamma$ to $\Sigma'$ that preserves $\{G_1, G_2\}$. Since $\Gamma$ is $2$PCI, there exists $n'\in N_{\S_{V(\Gamma)}}(R(G))$ such that $\Gamma^{n'}=\Sigma'$. It follows that $\Sigma=(\Sigma')^{(n')^{-1}n}$. Clearly, $(n')^{-1}n\in N_{\S_{V(\Gamma)}}(R(G))$, so $\Sigma$ is also $2$PCI.

For part (iv), letting $T=g^{-1}S$ for some $g\in S$, we have $\BCay(G,T)\cong\BCay(G,S)^{L_2(g)}$, where $L_2(g)\in N_{\S_{V(\Gamma)}}(R(G))$ is defined as $x_1\mapsto x_1, x_2\mapsto (g^{-1}x)_2, \forall x\in G$. By part (iii), $\BCay(G, T)$ is also $2$PCI. Clearly, $1_G\in g^{-1}S$.
\hfill\qed

\medskip
By Lemma~\ref{basic}~(i)-(ii), it is easy to see that $\BCay(G, S) \cong \BCay(G, T)$ if and only if there is an isomorphism from $\BCay(G, S)$ to $\BCay(G, T)$ keeping $\{G_1,G_2\}$ invariant (also see \cite[Lemma~4.4]{Arez}). Then Proposition~\ref{mCImPCI}~(2) implies the following corollary.

\begin{cor}\label{2pCI}
A $2$-PCayley (di)graph $\BCay(G, S)$ is $2$PCI if and only if for any $2$-PCayley (di)graph $\BCay(G, T)$ isomorphic to  $\BCay(G, S)$, there exist $\alpha\in\Aut(G)$ and $g\in G$ such that $T=g^{-1}S^{\alpha}$ or $T=(S^{-1})^{\alpha}g$, and is K$2$PCI if and only if for any $2$-PCayley (di)graph $\BCay(G, T)$ isomorphic to  $\BCay(G, S)$, there exist $\alpha\in\Aut(G)$ and $g\in G$ such that $T=g^{-1}S^{\alpha}$.
\end{cor}

\medskip
By definition, every K$2$PCI-graph is a $2$PCI-graph. However, the converse is not true. The following lemma provides some equivalent conditions for a $2$PCI-graph being K$2$PCI.

\begin{lem}\label{BCI2PCI}
Let $\Gamma=\BCay(G, S)$ be a $2$-PCayley graph of a group $G$. If $\Gamma$ is $2$PCI, then the following statements are equivalent.
\begin{enumerate}[\rm (1)]
  \item $\Gamma$ is a K$2$PCI-graph.
  \item $N_{\mathrm{Aut}(\Gamma)}(R(G))$ is transitive on $V(\Gamma)$.
  \item There exist $\alpha\in\Aut(G)$ and $g\in G$ such that $S^\alpha=S^{-1}g$.
\end{enumerate}
\end{lem}

\demo Assume that $\Gamma=\BCay(G,S)$ is $2$PCI. We first prove the equivalence between (1) and (2).
Suppose that $\Gamma$ is K$2$PCI. By Theorem~\ref{PCI-BabaiSimilar}, $N_{\mathrm{Aut}(\Gamma)}(R(G))$ is transitive on the bipartition $\{G_1, G_2\}$ of $V(\Gamma)$. Since $R(G)$ acts regularly on each part $G_i$, it follows that $N_{\mathrm{Aut}(\Gamma)}(R(G))$ is transitive on $V(\Gamma)$.
Conversely, assume that $N_{\mathrm{Aut}(\Gamma)}(R(G))$ is transitive on $V(\Gamma)$. Since $\Gamma$ is $2$PCI, by \cite[Theorem~4.5]{ZFYZ}, every semiregular subgroup of $\Aut(\Gamma)$ that is isomorphic to $G$ and has the orbit set $\{G_1, G_2\}$ on $V(\Gamma)$ is conjugate to $R(G)$ within $\Aut(\Gamma)$. Since $N_{\mathrm{Aut}(\Gamma)}(R(G))$ is transitive on $V(\Gamma)$, the sufficiency part of Theorem~\ref{PCI-BabaiSimilar} yields the conclusion that $\Gamma$ is K2PCI.

Next we shall prove that (1) and (3) are equivalent.
We first prove that (1) implies (3). By \cite[Proposition~2.1~(4)]{C-Z-F-Z}, we have $\Gamma\cong\BCay(G, S^{-1})$.
Since $\Gamma$ is K$2$PCI, by Corollary~\ref{2pCI}, there exist $\alpha\in\Aut(G)$ and $g\in G$ such that $S^{-1}=S^{\alpha}g^{-1}$. This proves (3). Now we prove that (3) implies (1). Since $\Gamma=\BCay(G,S)$ is $2$PCI, for any $2$-PCayley graph $\BCay(G, T)$ isomorphic to $\Gamma$, by Corollary~\ref{2pCI} there exist $\beta\in\Aut(G)$ and $h\in G$ such that either $T=h^{-1}S^{\beta}$ or $T=(S^{-1})^{\beta}h$. When $T=h^{-1}S^{\beta}$, Corollary~\ref{2pCI} gives that $\Gamma$ is K$2$PCI.
Now, consider the case when $T=(S^{-1})^{\beta}h$. By our assumption, $S^{\alpha}=S^{-1}g$ for some $\alpha\in\Aut(G)$ and $g\in G$. Then, $T=S^{\alpha\beta}g^{-\beta}h$.
Take $\gamma=\alpha\beta$ and $k=g^{-\beta}h$. Then $T=S^{\gamma}k$ for some $\gamma\in\Aut(G)$ and $k\in G$, and we obtain the desired result from Corollary~\ref{2pCI}.
\hfill\qed

\medskip
Let $G$ be an abelian group, and let $\Gamma=\BCay(G, S)$ be a $2$-PCayley graph of $G$. It is easy to see that $\a: g_1\mapsto g^{-1}_2$ and $g_2\mapsto g^{-1}_1$ for any $g\in G$, is an automorphism of $\Gamma$. Since $G$ is abelian, $g\mapsto g^{-1}$ is an automorphism of $G$. By Proposition~\ref{c-kernel}, $\a$ normalizes $R(G)$, and hence $N_{\mathrm{Aut}(\Gamma)}(R(G))$ is transitive on $V(\Gamma)$.

\begin{cor}\label{abelain2pCI}
Let $G$ be abelian. Then a $2$-PCayley graph $\BCay(G, S)$ is $2$PCI if and only if it is K$2$PCI.
\end{cor}


\subsection{Necessary and sufficient conditions for $2$PCI-groups}

In the following lemma, we give a necessary condition as well as a sufficient condition for $2$PCI-groups.
A finite group $G$ is called \emph{homogeneous} if every isomorphism between isomorphic subgroups of $G$ can be extended to an automorphism of $G$ (see~\cite{CherlinF}).

\begin{lem}\label{CompoBCI}
Let $G$ be a finite group. If $G$ is a $2$PCI-group, then for any $S\subseteq G$, $\BCay(\langle S\rangle, S)$ is $2$PCI, and any two subgroups of $G$ of the same order are equivalent under $\Aut(G)$.

Conversely, assume that for any $S\subseteq G$, $\BCay(\langle S\rangle, S)$ is $2$PCI and that any two subgroups of $G$ of the same order are equivalent under $\Aut(G)$. If $G$ is homogeneous, then $G$ is $2$PCI.
\end{lem}

\demo First assume that $G$ is $2$PCI. Let $S\subseteq G$ with $1\in S$ and let $H=\langle S\rangle$. Let $\BCay(H, T)$ be another $2$-PCayley graph on $H$ that isomorphic to $\BCay(H,S)$. By Lemma~\ref{basic}(iv), we may assume $1\in T$, and by Lemma~\ref{basic}(ii) we have $\BCay(G, S)\cong\BCay(G, T)$. Since $G$ is a $2$PCI-group, by Corollary~\ref{2pCI} we have $T=g^{-1}S^{\alpha}$ or $T=(S^{-1})^{\alpha}g$ for some $g\in G$ and $\alpha\in \Aut(G)$. As $1\in S$, it follows that either $g^{-1} \in T$ or $g\in T$, and by the connectedness of $\BCay(H, T)$, we have $H=\lg T\rg=\lg gT\rg=\lg Tg^{-1}\rg$. Thus, $H=H^\alpha$
because either $H^\alpha=\lg S^\alpha\rg=\lg gT\rg$ or $H^\alpha=\lg (S^{-1})^\alpha\rg =\lg Tg^{-1}\rg$.
Then the restriction of $\a$ on $H$ is an automorphism of $H$. Since $g\in H$,  Corollary~\ref{2pCI} implies that $\BCay(H, S)$ is $2$PCI, proving the first part of necessary condition.

Now let $M$ and $N$ be two subgroups of $G$ of the same order. Then $\BCay(M, M)\cong\BCay(N, N)\cong K_{|M|,|M|}$, and then $\BCay(G, M)\cong\BCay(G, N)$ by Lemma~\ref{basic}(ii). Since $G$ is a $2$PCI-group, by~Corollary~\ref{2pCI}, we have $M=g^{-1}N^{\alpha}$ or $M=(N^{-1})^{\alpha}g$ for some $g\in G$ and $\alpha\in \Aut(G)$. Since $M,N\leq G$, it follows that $M=N^\a$, completing the proof of necessary condition.

Conversely, assume that for any $S\subseteq G$ with $1\in S$, $\BCay(\langle S\rangle, S)$ is $2$PCI, and that
any two subgroups of $G$ of the same order are equivalent under $\Aut(G)$. Assume further that $G$ is homogeneous. We shall prove that $G$ is $2$PCI. To do this, let $\BCay(G, S)$ be a $2$-PCayley graph of $G$. Assume that $\BCay(G, T)$ is a $2$-PCayley graph of $G$ that is isomorphic to $\BCay(G,S)$. To finish the proof, by Corollary~\ref{2pCI} it suffices to show that there is $g\in G$ and $\alpha \in \Aut(G)$ such that $T=g^{-1}S^\alpha$ or $T=(S^{-1})^\alpha g$.
By Lemma~\ref{basic}~(i)-(ii), $\BCay(\langle T\rangle, T)$ and $\BCay(\langle S\rangle, S)$ are two isomorphic connected $2$-PCayley graphs.  Then $\langle T\rangle$ and $\langle S\rangle$ are two subgroups of $G$ with the same order. By the assumption, $\langle T\rangle$ and $\langle S\rangle$ are equivalent under $\Aut(G)$, that is, there is some $\sigma\in \Aut(G)$ such that $\langle S\rangle=\langle T\rangle^{\sigma}$. Then $\BCay(\langle S\rangle, T^{\sigma})$ is isomorphic to $\BCay(\langle T\rangle, T)$ under $\sigma^{-1}$, and since $\BCay(\langle T\rangle, T)\cong\BCay(\langle S\rangle, S)$, we have $\BCay(\langle S\rangle, T^{\sigma})\cong \BCay(\langle S\rangle, S)$. Since $\BCay(\langle S\rangle, S)$ is $2$PCI, there exist $g\in \langle S\rangle$ and $\alpha \in \Aut(\langle S\rangle)$ such that $T^{\sigma}=g^{-1}S^\alpha$ or $T^{\sigma}=(S^{-1})^\alpha g$. Thus $\beta:=\sigma \alpha^{-1}$ is an isomorphism from $\langle T\rangle$ to $\langle S\rangle$ such that $T^\beta=(T^\sigma)^{\alpha^{-1}}=(g^{-1})^{\alpha^{-1}}S$ or $T^\beta=(T^\sigma)^{\alpha^{-1}}=S^{-1}g^{\alpha^{-1}}$. As $G$ is homogeneous, there exists $\rho\in\Aut(G)$ such that $\beta=\rho|_{\langle T\rangle}$, the restriction of $\rho$ to $\langle T\rangle$. Therefore, $T^\rho=T^\beta=(g^{-1})^{\alpha^{-1}}S$ or $S^{-1}g^{\alpha^{-1}}$, and so $\BCay(G,S)$ is a $2$PCI-graph of $G$. This shows that $G$ is $2$PCI.\hfill\qed
\medskip





If $G$ is a direct product of $k$ cyclic groups of the same order then $G$ is said to be \emph{homocyclic} of rank $k$.

\begin{cor}\label{CompoBCIA}
A finite abelian group $G$ is $2$PCI if and only if every Sylow $p$-subgroup of $G$ is either elementary abelian or cyclic and
for any $S\subseteq G$, $\BCay(\langle S\rangle, S)$ is $2$PCI.
\end{cor}

\demo By \cite{Li6}, an abelian group is homogeneous if and only if every Sylow $p$-subgroup is homocyclic. So sufficiency is obvious by Lemma~\ref{CompoBCI}. For the necessity, if $G$ is an abelian $2$PCI-group, then by Lemma~\ref{CompoBCI}, it suffices to prove that every Sylow $p$-subgroup of $G$ is either elementary abelian or cyclic.
Actually, since $G$ is an abelian $2$PCI-group, by Lemma~\ref{CompoBCI}, for every Sylow $p$-subgroup $P$ of $G$, all subgroups of $P$ of order $p^2$ are isomorphic. It follows that every Sylow $p$-subgroup of $G$ is either elementary abelian or cyclic.\hfill\qed

\medskip 


By \cite[Theorem 5.5]{ZFYZ}, for $m\geq3$, every subgroup of a finite $m$PCI-group is $m$PCI. This, however is not true for $m=2$. Based on Lemma~\ref{CompoBCI}, in the next lemma, we shall present a necessary and sufficient condition for a subgroup of a finite $2$PCI-group being $2$PCI.\medskip

\begin{lem}\label{subBCI}
Let $G$ be a $2$PCI-group and $H\leq G$. Then $H$ is $2$PCI if and only if any two index $2$ subgroups of $H$ are equivalent under $\Aut(H)$.
\end{lem}

\demo Assume $G$ is $2$PCI. The necessity follows from Lemma~\ref{CompoBCI}.

For the sufficiency, let $\BCay(H,S)$ be a $2$-PCayley graph of $H$. Assume first that $\BCay(H,S)\ncong2K_{|H|/2,|H|/2}$, where $|H|$ is even. Then by~\cite[Lemma 4.6]{ZFYZ}, either $\BCay(H, S)$ or its complete bipartite complement is connected. Furthermore, by \cite[Corollary~4.7]{ZFYZ}, the complete bipartite complement of a $2$PCI-graph of $H$ is also a $2$PCI-graph of $H$. So we may assume that $\BCay(H,S)$ is connected. Then by Lemma~\ref{CompoBCI}, $\BCay(H,S)$ is $2$PCI.

Now assume that $\BCay(H,S)\cong2K_{|H|/2,|H|/2}$, where $|H|$ is even. By Lemma~\ref{basic}, we may assume $1\in S$. We claim $\BCay(H,S)\cong2K_{|H|/2,|H|/2}$ if and only if $S$ is an index $2$ subgroup of $H$.
If $S$ is an index $2$ subgroup of $H$, then $H=S\cup hS$ for some $h\in H\backslash S$, and $\BCay(S,S)\cong K_{|S|,|S|}\cong K_{|H|/2,|H|/2}$. It then follows that $\BCay(H,S)\cong2K_{|H|/2,|H|/2}$.
Conversely, if $\BCay(H,S)\cong2K_{|H|/2,|H|/2}$, then $|S|=|H|/2$. Notice that $1\in S$. Since $\BCay(H, S)$ is disconnected, it implies $\langle S\rangle\neq H$ and so $|S|\leq|\langle S\rangle|\leq |H|/2$. This implies that $S=\langle S\rangle$ is an index $2$ subgroup of $H$. The claim follows.

Let $\BCay(H,T)$ be another $2$-PCayley graph of $H$ isomorphic to $\BCay(H,S)$. By the claim above, $S$ and $T$ are index $2$ subgroups of $H$. Then, by the assumption, $S$ and $T$ are equivalent under $\Aut(H)$, and then there exists $\alpha\in \Aut(H)$ such that $T=S^\alpha$. This implies that $\BCay(H,S)$ is $2$PCI by~Corollary~\ref{2pCI}, completing the proof.\hfill\qed

\medskip

\begin{cor}\label{cor:charc-sub}
Let $G$ be a $2$PCI-group. If $H$ is a characteristic subgroup of $G$, then $H$ is $2$PCI.
\end{cor}

\demo Let $H$ be a characteristic subgroup of $G$ which contains index $2$-subgroups. For any two index $2$ subgroups $M$ and $N$ of $H$, since $G$ is $2$PCI, Lemma~\ref{CompoBCI} implies that there exists $\alpha\in \Aut(G)$ such that $N=M^\alpha$. As $H$ is characteristic in $G$, $\alpha|_{H}$, the restriction of $\alpha$ on $H$, is an automorphism of $H$. It follows that $M$ and $N$ are equivalent under $\Aut(H)$. Then by Lemma~\ref{subBCI}, $H$ is $2$PCI. \hfill\qed

\begin{lem}\label{chara2PCI}
Let $H$ be a characteristic subgroup of a $2$PCI-group $G$. Then $G/H$ is $2$PCI.
\end{lem}

\demo Assume $G$ is $2$PCI and $H$ is a characteristic subgroup of $G$. Let $\BCay(G/H, S)$ be a $2$-PCayley graph of $G/H$, and let $\BCay(G/H,T)$ be another $2$-PCayley graph of $G/H$ which is isomorphic to $\BCay(G/H,S)$. Let $\pi:G\rightarrow G/H$ denote the natural projection, which is a group homomorphism. Set $S_1:=\pi^{-1}(S)=\{g\mid g\in G, Hg\in S\}$. Clearly, if $g\in S_1$ then $hg\in S_1$ for any $h\in H$. Similarly, set $T_1:=\pi^{-1}(T)=\{g \mid g\in G, Hg\in T\}$. We claim that $\BCay(G,S_1)$ is isomorphic to the lexicographic product $\BCay(G/H,S)[nK_1]$, where $n = |H|$.

Consider an edge $\{(Hg)_1, (Hsg)_2\}$ in $\BCay(G/H,S)$, where $Hs\in S$. Then  $hs\in S_1$ for every $h\in H$. For any $h, h'\in H$. We will show $\{(hg)_1, (h'sg)_2\}$ is an edge of $\BCay(G,S_1)$.
Since $H$ is characteristic in $G$, there exists $h_1\in H$ such that $h_1s = sh$, or equivalently, $h_1^{-1}s = sh^{-1}$. Then, we have $(h'h_1^{-1}shh^{-1})hg = h'sg$, where $h'h_1^{-1}shh^{-1}=h'h_1^{-1}s\in S_1$.
Consequently, $\{(hg)_1, (h'sg)_2\}\in E(\BCay(G,S_1))$. This shows that the subgraph of $\BCay(G,S_1)$ induced by $(Hg)_1\cup (Hsg)_2$ is isomorphic to the complete bipartite graph $K_{n,n}$. As $\BCay(G,S_1)$ is bipartite, for every $g\in G$ and $1\leq i\leq 2$, the set $(Hg)_i$ induces an empty graph in $\BCay(G,S_1)$. It follows that $\BCay(G,S_1)\cong\BCay(G/H,S)[nK_1]$.

Since $G$ is a $2$PCI-group,  by~Corollary~\ref{2pCI}, there exist $\alpha\in\Aut(G)$ and $g\in G$ such that either $T_1=g^{-1}S_{1}^{\alpha}$ or $T_1=(S_{1}^{-1})^{\alpha}g$. Also as $H$ is a characteristic subgroup of $G$, let $\overline{\alpha}:=\alpha|_H$, the restriction of $\alpha$ on $H$. Then $\overline{\alpha}$ is an automorphism of $G/H$ satisfying $(Hg)^{\overline{\alpha}}=Hg^{\alpha}$ for every $g\in G$. Now we have $T=T_1^{\pi}=(g^{-1})^{\pi}(S_1^{\alpha})^{\pi}=(g^{-1})^{\pi}S^{\overline{\alpha}}$ or $T=T_1^{\pi}=(S_1^{-1})^{\alpha\pi}g^{\pi}=(S^{-1})^{\overline{\alpha}}g^{\pi}$. By Corollary~\ref{2pCI}, $\BCay(G/H,S)$ is $2$PCI. This completes the proof.
\hfill\qed

\subsection{Proof of Theorem~\ref{Sylow2pci}}

In this section, we shall prove Theorem~\ref{Sylow2pci}. Before proceeding, we introduce two results.
A group $G$ is said to be an \emph{iso-group} if any two subgroups of $G$ of the same order are isomorphic. The following proposition gives a complete description of finite iso-groups, which is due to Zhang~\cite{Z}.


\begin{prop}[{\cite[Theorem 8]{Z}}]\label{iso-group}
If $G$ is an iso-group then $G$ is isomorphic to one of the following groups:
\begin{enumerate}[\rm(1)]
  \item $K\rtimes H$; $(|H|,|K|)=1$, $K$ is a normal nilpotent subgroup and each Sylow subgroup of $K$ is cyclic, elementary, quaternion or extra-special of order $p^3$ and of exponent $p$; all Sylow subgroups of $H$ are cyclic.
  \item $K\rtimes((Q_8 \times H):\langle\alpha\rangle)$; $H$ and $K$ are as above, $H$, $K$, $Q_8$ and $\langle\alpha\rangle$ are of pairwise coprime order, $\langle\alpha\rangle$ is a Sylow $3$-subgroup of $G$ acting nontrivially on $Q_8$ if $\alpha \neq 1$.
  \item $K\rtimes(\mathrm{SL}(2,5) \times H)$; $H$, $K$ and $\mathrm{SL}(2,5)$ are of pairwise coprime order and $H$ and $K$ are as in~\rm(1).
  \item $(\mathrm{SL}(2,2^n) \times L)\rtimes\langle d\rangle$; $n \geq 2$, $\mathrm{SL}(2,2^n)$ and $L\langle d\rangle$ are of coprime order and $L\langle d\rangle$ is of type~\rm(1).
\end{enumerate}
\end{prop}

A group is called an {\em FIF-group} if for any two elements $x, y$ of $G$ of the same order, there exists $\a\in\Aut(G)$ such that $x^\a=y$ or $x^\a=y^{-1}$.

\begin{lem}\label{lem:FIF}
Every $2$PCI-group is an FIF-group.
\end{lem}

\demo Let $G$ be a $2$PCI-group. Let $s, t$ be two elements of $G$ of the same order. Let $S=\{1, s\}$ and $T=\{1, t\}$. Then both $\BCay(G, S)$ and $\BCay(G, T)$ are isomorphic to a union of cycles of length $2|s|=2|t|$. Since $G$ is a $2$PCI-group, by Corollary~\ref{2pCI}, there exist $g\in G$ and $\alpha\in \Aut(G)$ such that
\[{\rm either}\ gS^{\alpha}=T\ {\rm or}\ (S^{-1})^\alpha g^{-1}=T.\]
For the former, we have $g=1$ or $t$: for $g=1$ we have $s^\a=t$ and for $g=t$ we have $s^\a=t^{-1}$.
For the latter, we have $g=1$ or $t^{-1}$, and we also have $s^\a=t$ or $t^{-1}$. This implies that $G$ is an FIF-group. \hfill\qed

\medskip
The following result gives a characterization of FIF-groups, which is due to Li and Praeger.

\begin{prop}[{\cite[Corollary~1.3]{Li2}}]\label{FIF-group}
If $G$ is a FIF-group. Then the following hold:
\begin{enumerate}[\rm(1)]
  \item If $G$ is a non-abelian simple group then $G$ is one of $\mathrm{PSL}_2(q)$ for $q=5,7,8$ or 9 , $\mathrm{PSL}_3(4)$, $\mathrm{Sz}(8), \mathrm{M}_{11}$ and $\mathrm{M}_{23}$.
  \item If $G$ is nilpotent then each Sylow subgroup of $G$ is either a homocyclic group, $\mathrm{Q}_8$, or belongs to a class of non-abelian $2$-group.
  \item If $G$ is soluble then $G=A: B$ with $(|A|,|B|)=1$, where $A$ is a nilpotent FIF-group and every Sylow subgroup of $B$ is cyclic or $\mathrm{Q}_8$.
  \item If $G$ is insoluble then $G=A \times B$, where $A, B$ have coprime orders, $A$ is a soluble FIF-group and $B$ either is one of the simple groups in (1), or is $\mathrm{SL}_2(q)$ for $q=5,7$ or 9 , or is $(C \times \mathrm{Sz}(8)): \mz_{3^sm}$, where $m, s \geq 1$.
\end{enumerate}
\end{prop}
\medskip

\noindent{\bf Proof of Theorem~\ref{Sylow2pci}:}\ Assume $G$ is a $2$PCI-group. 
By Lemma~\ref{CompoBCI}, any two subgroups of $G$ of the same order are equivalent under $\Aut(G)$. It follows that $G$ is an iso-group, and so $G$ is one of the groups given in Proposition~\ref{iso-group}. We claim that $G$ is solvable. Suppose on the contrary that $G$ is nonsolvable. Then $G$ is isomorphic to one of the groups given in Proposition~\ref{iso-group}(3) or (4).

If Proposition~\ref{iso-group}(3) happens, then we may let $G=K\rtimes(\mathrm{SL}(2,5) \times H)$, where $H$, $K$ and $\mathrm{SL}(2,5)$ are of pairwise coprime order and $H$ and $K$ are as in Proposition~\ref{iso-group}(1). Clearly, $\mathrm{SL}(2,5)$ has no index $2$ subgroup, and by Lemma~\ref{subBCI}, $\mathrm{SL}(2,5)$ is a $2$PCI-group. Since  $Z(\mathrm{SL}(2,5))$ is characteristic in $\mathrm{SL}(2,5)$, by Lemma~\ref{chara2PCI}, $\mathrm{PSL}(2,5)=\mathrm{SL}(2,5)/Z(\mathrm{PSL}(2,5))$ is a $2$PCI-group.
Note that $\mathrm{PSL}(2,5)\cong \mathrm{A}_5$. Take $S=\{(13245),(15423),(24)(35),(12)(45),(15)(23),\mathrm{Id}(G)\}$ and let $\Gamma=\BCay(\mathrm{A}_5, S)$. By {\sc Magma}~\cite{magma}, $\Aut(\Gamma)$ has a semiregular subgroup which is isomorphic to $R(\mathrm{A}_5)$ and has the same orbits with $R(\mathrm{A}_5)$, but is not conjugate to $R(\mathrm{A}_5)$ in $\Aut(\Gamma)$. By Theorem~\cite[Theorem~4.5]{ZFYZ}, $\Gamma$ is not $2$PCI, and hence $\mathrm{PSL}(2,5)$ is not a $2$PCI-group, a contradiction.

If Proposition~\ref{iso-group}~(4) happens, then we may let $G=(\mathrm{SL}(2,2^n) \times L)\rtimes\langle d\rangle$; $n\geq 2$, $\mathrm{SL}(2,2^n)$ and $L\langle d\rangle$ are of coprime order and $L\langle d\rangle$ is of type~\rm(1) of Proposition~\ref{iso-group}.
This implies that $\mathrm{SL}(2,2^n)$ is the only nonsolvable chief factor. Recall that $G$ is an FIF-group. From~\cite[Corollary 1.3~(1) and (4)]{Li2} we must have $G=A\times B$, where $A=\mathrm{SL}(2, 8)$, $(|A|, |B|)=1$, and $B$ is a solvable FIF-group. So $A=\mathrm{SL}(2, 8)$ is characteristic in $G$. By Corollary~\ref{cor:charc-sub}, $\mathrm{SL}(2, 8)$ is $2$PCI. Note that $\mathrm{SL}(2, 8)$ has a maximal subgroup which is isomorphic to the Frobenius group $F_8=\Z_2^3: \Z_7$. Since $F_8$ has no index $2$ subgroups, by Lemma~\ref{subBCI}, $F_8$ is $2$PCI.

Let $M=\langle x,y,z,t\rangle$, where $x=(2687453)$, $y=(13)(24)(57)(68)$, $z=(12)(34)(56)(78)$, and $t=(15)(26)(37)(48)$ in $\mathrm{S}_8$. Then $M\cong F_8$. Take $S=\{(1327846),(18)(27)(36)(45),(1742365),\mathrm{Id}(K)\}$ and let $\Gamma=\BCay(K,S)$. By {\sc Magma}~\cite{magma}, $\Aut(\Gamma)$ has a semiregular subgroup that is isomorphic to $R(M)$ and has the same orbits with $R(M)$, but it is not conjugate to $R(K)$ in $\Aut(\Gamma)$. By Theorem~\cite[Theorem~4.5]{ZFYZ}, $\Gamma$ is not a $2$PCI-graph, and thus $F_8$ is not $2$PCI, a contradiction.

By now, we have shown that $G$ is solvable, and so $G$ is isomorphic to one of the groups given in Proposition~\ref{iso-group}(1) or (2). Let $P$ be a Sylow $p$-subgroup of $G$. Inspecting the groups in Proposition~\ref{iso-group}(1) or (2), we may conclude that either $P$ is elementary abelian, or one of the following happens:
\begin{itemize}
  \item $P$ is cyclic,
  \item $P\cong Q_8$,
  \item $P\cong SP:=\langle a,b,c\ |\ a^p=b^p=c^p=1,[a,b]=c,[a,c]=[b,c]=1\rangle (p\geq3)$.
\end{itemize}
If $p>2$, then $P$ has no index $2$-subgroups, and if $p=2$, then any subgroup $Q$ of $P$ either is isomorphic to $Q_8$, or is elementary abelian or cyclic, and so any two index $2$ subgroups of $Q$ are equivalent under $\Aut(Q)$. By  Lemma~\ref{subBCI}, $P$ is $2$PCI, and
by Lemma~\ref{lem:FIF}, $P$ is an FIF-group. By Corollary~\ref{FIF-group}~(2) and (3), $P\ncong SP$ and hence $P$ is elementary abelian, or cyclic, or $Q_8$.

First assume $p=2$. Let $K=\langle x\rangle\cong\mz_8$ and $S=\{1_K, x, x^2,x^5\}$.
Let $\Gamma=\BCay(K,S)$. By {\sc Magma}~\cite{magma}, $R(K)\unlhd \Aut(\Gamma)$, and there is another semiregular subgroup of $\Aut(\Gamma)$ that is isomorphic to $K$ and has same orbits with $R(K)$, but it is not conjugate to $R(K)$ in $\Aut(\Gamma)$. By Theorem~\cite[Theorem~4.5]{ZFYZ}, $\Gamma$ is not a $2$PCI-graph, and thus $\mz_8$ is not $2$PCI.
It follows that $\Z_{2^n}$ is not $2$PCI-group for any $n\geq 3$. Thus, $P$ is elementary abelian, $\Z_4$, or $Q_8$.

Secondly assume $p=3$. Let $M=\langle x\rangle\cong\mz_{3^3}$ and let
$S=\{1_M, x^{11}, x^{-7}, x^2, x^9, x^{-9}, x^4, x^{-2}\}$.
Let $\Gamma=\BCay(M, S)$. By {\sc Magma}~\cite{magma}, $\Aut(\Gamma)$ has a semiregular subgroup that is isomorphic to $R(M)$ and has the same orbits with $R(M)$, but not conjugate to $R(M)$ in $\Aut(\Gamma)$. By Theorem~\cite[Theorem~4.5]{ZFYZ}, $\Gamma$ is not a $2$PCI-graph, and thus $\mz_{3^3}$ is not $2$PCI. By Corollary~\ref{abelain2pCI}, $\Z_{3^n}$ is not $2$PCI-group for any $n\geq 3$. Let $R=\langle a\rangle\times\langle b\rangle\times\langle c\rangle\cong\Z_3^3$ and $S=\{1_R, a,b,c,a^2,b^2,ac^2,b^2c^2,a^2bc\}$. By {\sc Magma}, $\BCay(R, S)$ is not a K$2$PCI-graph, and hence $\Z_3^n$ is not $2$PCI for all $n\geq 3$. It follows that $P\cong\Z_3, \Z_9$ or $\Z_3^2$.

Finally assume $p>3$. By~\cite[Proposition 4.10]{Arez}, $\Z_{p^2}$ is not a K$2$PCI-group, and by Corollary~\ref{abelain2pCI}, $\Z_{p^2}$ is not $2$PCI-group. Thus, $P$ is elementary abelian.

Now we know that $P$ is isomorphic to one of the groups listed in Theorem~\ref{Sylow2pci}.\hfill\qed





\subsection{Proof of Theorem~\ref{BCIEA}}

We first introduce the following proposition which is due to Morris and Spiga.

\begin{prop}[{\cite[Theorem 1.1]{Morris}}]\label{BCI}
Let $G$ be a finite group. Then one of the following holds
\begin{itemize}
  \item [\rm(1)] $G$ admits a $2$-PCayley graphical representation;
  \item [\rm(2)] $G$ is an abelian group and there exists $S \subseteq G$ with $R(G) \rtimes\langle\iota\rangle=\Aut(\BCay(G, S))$, where $\iota$ is the map $g_1\mapsto g^{-1}_2$ and $g_2\mapsto g_1^{-1}$ for any $g\in G$;
  \item [\rm(3)] $G$ is one of the twenty-two exceptional groups appearing in Table \rm1.
\end{itemize}
\end{prop}
\begin{table}[h]
\begin{center}
\begin{tabular}{c|c|c|c}
\hline Line&Order & Group &K$2$PCI\\
\hline 1&3 & $\Z_3$& Y \\
2&4 & $\Z_2^2$ or $\Z_4$ & Y\\
3&5 & $\Z_5$& Y \\
4&6 & $\Z_6$ or $D_6$& Y \\
5&7 & $\Z_7$ & Y\\
6&8 & $\Z_2^3$ or $Q_8$ or & Y\\
&&$\Z_4 \times \Z_2$ or $D_8$ & N\\
7&9 & $\Z_3^2$ & Y\\
8&10 & $D_{10}$&Y \\
9&12 & Alt$(4)$ or $D_{12}$ or $\langle x, y \mid x^6=y^4=1, x^3=y^2, y^{-1} x y=x^{-1}\rangle$&N \\
10&14 & $D_{14}$&N \\
11&16 & $\Z_2^4$ or &Y\\
&&$\Z_4 \times \Z_2^2$ or $Q_8 \times \Z_2$&N \\
12&18 &$\langle e_1, e_2, x \mid e_1^3=e_2^3=x^2=[e_1, e_2]=1, e_1^x=e_1^{-1}, e_2^x=e_2^{-1}\rangle$&N \\
13&32 & $\Z_2^5$&Y\\
\hline
\end{tabular}
\caption{Exceptional examples in Proposition~\ref{BCI}.}\label{tb:bci}
\end{center}
\end{table}

\noindent{\bf Proof of Theorem~\ref{BCIEA}:}\ (1)\ If $G$ is abelian, then by Corollary~\ref{abelain2pCI}, Corollary~\ref{CompoBCIA} and Theorem~\ref{Sylow2pci}, we have $G$ is K$2$PCI (namely, BCI) if and only if every Sylow $3$-subgroup of $G$ is isomorphic to $\Z_3, \Z_3\times\Z_3$ or $\Z_9$, every Sylow $p$-subgroup ($p\neq 3$) of $G$ is either elementary abelian, or is isomorphic to $\Z_4$, and for any $S\subseteq G$, $\BCay(\langle S\rangle, S)$ is K$2$PCI.

(2)\ Suppose that $G$ is a non-abelian K$2$PCI group. Then $G$ is in one of the two classes of groups in Proposition~\ref{BCI}(1) and (3). Since $G$ is a non-abelian K$2$PCI group, by Lemma~\ref{BCI2PCI}, every $2$-PCayley graph of $G$ is transitive. This implies that Proposition~\ref{BCI}(1) cannot happen. So $G$ is one of the non-abelian groups in Table~\ref{tb:bci}. Note that $G$ is also a $2$PCI-group since it is a non-abelian K$2$PCI-group. By Theorem~\ref{Sylow2pci}, we have $G$ is not isomorphic to one of the following groups:
\[D_8, Q_8 \times \Z_2, \Z_4 \times \Z_2^2, \Z_4 \times \Z_2.\]

By \cite[Remark~1 \& Proposition~11]{EP}, if $n>5$, then $D_{2n}$ has a $2$-PCayley graph which is not vertex-transitive, and so by Lemma~\ref{BCI2PCI}, $D_{2n}$ is not a K$2$PCI group for each $n>5$.  So if $G$ is dihedral, then $G\cong D_6$ or $D_{10}$.

Now we let
\[\begin{array}{l}
G_1:=\mathrm{A}_4, S_1:=\{(143),(234), (13)(24),\mathrm{Id}(G)\},\\
G_2:=\langle x, y \mid x^6=y^4=1, x^3=y^2, y^{-1} x y=x^{-1}\rangle, S_2:=\{x^2,xy, xy^{-1}, x^{-1}, \mathrm{Id}(G), yx\},\\
G_3:=\langle e_1, e_2, x \mid e_1^3=e_2^3=x^2=[e_1, e_2]=1, e_1^x=e_1^{-1}, e_2^x=e_2^{-1}\rangle, S_3:=\{ e_1, e_2, x, e_1x, e_{1}^{-1}, xe_1e_2\}.
\end{array}\]
With help of the software {\sc Magma}~\cite{magma}, all the $2$-PCayley graphs $\BCay(G_i, S_i)$ for $i=1,2,3$ are not vertex-transitive, and hence all these three groups $G_1, G_2, G_3$ are not K$2$PCI-groups by Theorem~\ref{PCI-BabaiSimilar}.

By now, we have shown that $G\cong D_6, D_{10}$ or $Q_8$. 
By Magma~\cite{magma} we know that $D_6$, $Q_8$ and $D_{10}$ are all K$2$PCI-group. This proves part (2).\hfill\qed

\medskip
We also find all K$2$PCI-groups among the groups in Table~\ref{tb:bci}. In the last column of Table~\ref{tb:bci}, `Y' means that the groups in the same line are K$2$PCI, while `N' means the groups in the same line are not K$2$PCI. We now explain how to obtain this. Actually, we can prove that $\Z_2^2$, $\Z_2^3$, $\Z_2^4$ and $\Z_2^5$ are all K$2$PCI-groups (see Lemma~\ref{Z24}). By~\cite[Corollary 4.9]{Arez}, we know that $\Z_3,\Z_5,\Z_7$ are all K$2$PCI-groups.
With the help of {\sc Magma}~\cite{magma}, we know that $\Z_4,\Z_6,\Z_3^2$ are also K$2$PCI.

\begin{lem}\label{Z24}
$\Z_2^n$ is K$2$PCI for each $0\leq n\leq 5$.
\end{lem}

\demo Since $\Z_4$ is K2PCI, Corollary~\ref{abelain2pCI} implies that  $\Z_4$ is 2PCI. By Lemma~\ref{CompoBCI},  $\Z_1$ and $\Z_2$ are 2PCI, and again by Corollary~\ref{abelain2pCI}, $\Z_1$ and $\Z_2$ are K2PCI. To prove our lemma, by Theorem~\ref{BCIEA}, it is enough to prove that for each $2\leq n\leq 5$, every connected $2$-PCayley graph of $\Z_2^n$ is K$2$PCI.

Assume $H=\Z_2^n$ for some $2\leq n\leq 5$.
Let $\Gamma=\BCay(H,S)$ be a connected $2$-PCayley graph of $H$. By ~\cite[Corollary~4.7]{ZFYZ} we know that $\Gamma$ is K$2$PCI if and only if its complete bipartite complement is K$2$PCI. We may assume that $1\leq|S|\leq |H|/2$.
By Lemma~\ref{basic}(iv), we may further assume the $1_H\in S$
. As $\Gamma$ is connected, by Lemma~\ref{basic}(i), we have $\langle S\rangle=H$, and hence $|S|\geq n+1$. This implies that $n\geq3$.

Suppose that $S\setminus\{1_H\}$ is a minimum generating subset of $H$. Let $\BCay(H, T)$ be another $2$-PCayley graph of $H$ that is isomorphic to $\Gamma$. We may assume $1_H\in T$. Then $T\setminus\{1_H\}$ is also a minimum generating subset of $H$. So there is an $\alpha\in\Aut(H)$ such that $T=S^{\alpha}$. It follows from Corollary~\ref{2pCI} that $\Gamma$ is K2PCI. With these arguments, in the following, we may assume that $S\setminus\{1_H\}$ is not a minimum generating subset of $H$. Then either $n=4$ and $6\leq |S|\leq 8$, or $n=5$ and $7\leq |S|\leq 16$.

\medskip
\noindent{\bf Case~1.}\ $n=4$ and $6\leq |S|\leq 8$.

Let $H=\langle a\rangle\times\langle b\rangle\times\langle c\rangle\times\langle d\rangle\cong\Z_2^4$. Due to $\langle S\rangle=H$, we may assume that $a,b,c,d\in S$. Let $T=\{1_H, a,b,c,d\}$. By {\sc Magma}~\cite{magma}, there exists $\a\in\Aut(H)$ such that $S^\a=T\cup T'$, where $T'\in\Omega$ with
\[
\begin{array}{ll}
\Omega=&\{\{ab\}, \{abc\}, \{abcd\}, \\
&\{ab,ac\}, \{ab,bcd\}, \{ab,cd\}, \{abc,abd\},\\
& \{ab,ac,ad\}, \{ab,ac,bc\}, \{ab,ac,bd\}, \{ab,ac,bcd\},\{abc,abd,acd\}\}.
\end{array}
\]
This implies that $\BCay(H, S)\cong \BCay(H, T\cup T')$, where $T'\in\Omega$. By {\sc Magma}~\cite{magma}, for every $T'\in\Omega$, every semiregular subgroup of $\Aut(\BCay(H, T\cup T'))$ that is isomorphic to $\Z_2^4$ and has the same orbits as $R(H)$ is conjugate to $R(H)$ in $\Aut(\BCay(H, T\cup T'))$. It follows that $\BCay(H, T\cup T')$ is K2PCI for every $T'\in \Omega$. Thus, every connected $2$-PCayley graph of $H$ of valency $6$, $7$ or $8$ is K2PCI.

\medskip
\noindent{\bf Case~2.}\ $n=5$ and $7\leq |S|\leq 16$.

In this case, $H\cong\Z_2^5$. For convenience, let $R(H)$ be a semiregular permutation group on the set $\Omega=\{1, 2,\ldots, 64\}$ with two orbits $\Omega_1=\{1, 2, \ldots, 32\}$ and $\Omega_2=\{33, 34, \ldots, 64\}$. Let $\Gamma=\BCay(H, S)$ be a graph $\mathcal{G}(O)$ with vertex set $\Omega$ and edge set $\{\{1, 32+i\}^{R(h)}\mid i\in O, h \in H\}$, where
\[O\in \Delta:=\{O\mid O\subseteq \{1,2,\ldots,32\},  7\leq |O|\leq 16, 1\in O, \text{$\mathcal{G}(O)$ is connected} \}.\]



Let $N=N_{\S_{\Omega}}(R(H))$ and let $K$ be the kernel of $N$ acting on $\{\Omega_1, \Omega_2\}$.
These groups $R(H)$, $N$ and $K$ can be easily constructed in {\sc Magma}~\cite{magma}.
Clearly,
\[
\Delta=\bigcup_{t=7}^{16}\Delta_t, \text{ where } \Delta_t:=\{O\mid  O  \in \Delta, |O|=t \}.
\]
Next we shall find a representative from each orbit of $K$ on $\Delta_t$ for $7\leq t\leq 16$, and all these representatives form a set, say $\Sigma$. This can be done with help of the software package {\sc Orbiter}~\cite{Orbiter}. Furthermore, for any two different $O_1, O_2\in\Sigma$, by {\sc Magma}~\cite{magma}, we can verify that $\mathcal{G}(O_1)\ncong \mathcal{G}(O_2)$. This implies that for every $\mathcal{G}(O')$ with $O' \in \Delta$ such that $\mathcal{G}(O')\cong \mathcal{G}(O)$, there exists $k\in K$ such that $\mathcal{G}(O')=\mathcal{G}(O)^k$.
By definition, we know that every connected $2$-PCayley graph $\BCay(H, S)$ of $H\cong\Z_2^5$ with $1\in S$ and $7\leq|S|\leq 16$ is K$2$PCI, completing the proof.  \qed

\section*{Acknowledgments}
This work was partially supported by the National Natural Science Foundation of China (12271024, 12331013, 12301461, 12425111).
	
	

{}

\end{document}